\documentclass[a4paper,12pt]{article}

\usepackage{amsmath,latexsym,amssymb}
\usepackage[german,english]{babel}
\usepackage{color} 
 
\newtheorem{thm}{Theorem}
\newtheorem{prop}{Proposition}[section]
\newtheorem{definition}{Definition}
\newtheorem{remark}{Remark}
\newtheorem{lemma}[prop]{Lemma}

\newtheorem{notation}{Notation}

\newcommand{\proof}[1][]{{\it Proof#1: }}
\newcommand{\qed}[1][3mm]{\hspace*{\fill} $\Box$ \vspace{#1}}

\newcommand{\vv}{{\mathbf v}}

\newcommand{\bN}{{\mathbb N}}
\newcommand{\bZ}{{\mathbb Z}}
\newcommand{\bC}{{\mathbb C}}
\newcommand{\bP}{{\mathbb P}}

\newcommand{\ul}{\underline}
\newcommand{\beq}{ \stackrel{B}{\sim}}
\newcommand{\aeq}{ \stackrel{A}{\sim}}
\newcommand{\baeq}{ \stackrel{BA}{\sim}}

\renewcommand{\a}{\alpha}
\renewcommand{\b}{\beta}
\newcommand{\cg}{\gamma}

\newcommand{\s}{\sigma}

\newcommand{\ra}{\rightarrow}

\newcommand{\inv}{^{^{-1}}}

\newcommand{\hfami}{{\cal H}}

\newcommand{\sC}{\ensuremath{\mathcal{C}}}
\newcommand{\sS}{\ensuremath{\mathcal{S}}}
\newcommand{\sT}{\ensuremath{\mathcal{T}}}
\newcommand{\sB}{\ensuremath{\mathcal{B}}}

\newcommand{\de}{\delta}

\newcommand{\Dn}{{\rm D}_n}
\newcommand{\Dm}{{\rm D}_m}

\newcommand{\cutoff}[1]{}

\begin{document}

\title{Irreducibility of the space of dihedral 
covers of the projective line of a given 
numerical type}

\author{Fabrizio Catanese, Michael L\"onne, Fabio Perroni
\thanks{The present work took place in the realm of the DFG
Forschergruppe 790 "Classification of algebraic
surfaces and compact complex manifolds". }
}

\date{\today}

\maketitle

\hfill {\em Dedicated to the memory of Giovanni Prodi.}

\bigskip

\begin{abstract}
We show in this paper that the set of irreducible 
components of the family of Galois  coverings of
$\bP^1_{\bC}$ with Galois group isomorphic to $\Dn$ 
is in bijection with the set of
possible numerical types.

In this special case the numerical type is the equivalence class
(for automorphisms of $\Dn$) of the function 
which to each conjugacy class $\sC$ in
$\Dn$ associates the number of branch points whose local monodromy lies
in the class $\sC$.

\end{abstract}

KEY WORDS:  Moduli spaces of curves, branched coverings of Riemann surfaces, Hurwitz equivalence, braid groups,
monodromy.

2000 MATHEMATICS SUBJECT CLASSIFICATION: 14H10, 20F36, 57M12.

%
%
%
%
%


\section*{Introduction}

The theory of covering spaces was invented to 
clarify the concept of an algebraic
function and its polydromy.

In the modern terminology one can describe an 
algebraic function $f$ on an algebraic curve $Y$ 
as a rational
function
$f$ on a projective curve $X$ admitting a non constant
morphism $p : X \ra Y$, and such that $f$ 
generates the field extension
$ \bC(Y) \subset \bC(X)$.

The easiest example would be the one where $Y 
=\bP^1 = \bP^1_{\bC}$ and $f= \sqrt P(x)$,
    $P$ being a square free polynomial.

$f$ is in general polydromic, i.e., many valued as a function on $Y$, and
going around a closed loop we do not return to the same value.
It is a theorem of Weierstrass that $f$ is a rational
function on $Y$ iff $f$ is monodromic, i.e., there is no polydromy.

For strange reasons (but remember the famous 
explanation  `Lucus a non lucendo', the grove 
has a similar
name to light because there is no light) what 
should be called polydromy is nowadays called 
monodromy.

Given $p : X \ra Y$ as above there is a finite 
set $\sB \subset Y$, called the branch locus,
such that, setting  $Y^* : = Y \setminus \sB, X^* 
: = p^{-1} (Y^* )$, then $p$ induces a covering
space $X^* \ra Y^*$ which is classified by its 
monodromy $\mu$, which is a homomorphism
of the fundamental group (Poincar\'e group) of 
$Y^*$, $\pi_1 (Y^*, y_0)$, into the group
of permutations of the fibre $  p^{-1} (y_0)$.

If $X$ is irreducible, and $\de$ is the degree of 
$p$, then the image of the monodromy is a 
transitive
subgroup of $\mathfrak S_{\de}$, and
conversely Riemann's existence theorem asserts that for any homomorphism
$ \mu \colon  \pi_1 (Y^*, y_0) \ra \mathfrak S_{\de}$ 
with transitive image we obtain a morphism $ p \colon X \ra Y$ as above inducing
the given monodromy $\mu$, hence also a corresponding  algebraic 
function
on $Y$ with branch set contained in $\sB$.

Riemann's existence theorem is a very powerful 
but not constructive result: it is similar in 
spirit
to the non constructive argument which shows that 
any $n \times n$ matrix $A$ satisfies
a polynomial equation of degree at most $n^2$; 
while the theorem of Hamilton Cayley
constructs such a polynomial equation of degree 
$n$, namely, the characteristic polynomial $P_A$ 
of $A$.
Although $P_A$ is not the polynomial of minimal 
degree which gives zero when evaluated
on $A$, it has the  advantage that it varies well 
with $A$ if $A$ varies in a family.

Similarly, one can consider families of algebraic 
functions, or, equivalently, families
of morphisms $X_t \ra Y_t$ of algebraic curves, 
and a natural question is whether
a given  parameter space $T$
is irreducible: for this type of question 
Riemann's existence theorem plays a crucial role.

Usually one splits the above question by 
considering families where the branch locus has a 
fixed
cardinality, obtaining in this way a 
stratification of the parameter space (the strata 
are often called Hurwitz
spaces, see \cite{Fulton}); and then asking  which  strata are irreducible.

The archetypal result is the theorem of 
L\"uroth-Clebsch and Hurwitz, showing that
simple coverings of the projective line form an 
irreducible variety (see \cite{Clebsch}, 
\cite{Hurwitz},
cf. also \cite{BC} for a simple proof).
Here, simple means that the local monodromies 
(image under $\mu$ of small loops around the 
branch
points) are transpositions. The theorem of 
L\"uroth-Clebsch has been extended to projective 
curves
$Y$ of higher genus by several authors 
(\cite{GHS}), and there are variants 
(\cite{Kluitmann},\cite{Bronek1},
\cite{Ka1}, \cite{Ka2}, 
\cite{Vetro1},\cite{Vetro2},\cite{Vetro3}) where 
for one
or two distinguished branch points the local 
monodromy can be chosen to be a different type of 
permutation,
or where one replaces the symmetric group by Weyl 
groups and the transpositions by reflections.

   Observe that one can factor the monodromy $ \mu 
\colon  \pi_1 (Y^*, y_0) \ra \mathfrak S_{\de}$
through a surjection onto a finite group $G$ 
followed by a permutation representation of $G$, 
i.e., an injective
homomorphism $ G \ra \mathfrak S_{\de}$ with transitive image.

Geometrically this amounts
to construct a morphism $Z \ra Y$ (the Galois 
closure of $p$) such that $G$ acts on $Z$ with 
quotient $ Z/G
\cong Y$, and such that $X$ is obtained as the 
quotient of $Z$ by a non normal subgroup $ H$ of 
$G$,
and we have the factorization
$$ Z \ra Z/ H = X \ra Z/G = Y. $$

In this way one separates the investigation of 
algebraic functions into two parts: the study of 
Galois covers
$ Z \ra Y$, and the study of intermediate covers.

  The study of Galois covers is however also of interest in
itself, since inside the moduli space $\mathfrak 
M_g$ of curves $X$ of genus $g \geq 2$ we have
the closed proper algebraic subset of curves 
having a nontrivial group of automorphisms,
and one would like to understand, given a finite 
group $G$, which are the irreducible components
of the algebraic subset $\mathfrak M_{g;G}$ of 
curves $X$ admitting $G$ as a subgroup
of their group of automorphisms.

The action of $G$ on the curve $X$ gives rise to 
a morphism $p \colon  X \ra X/G = Y$,
and the geometry of $p$ encodes several discrete 
invariants which distinguish the
  irreducible components of $\mathfrak M_{g;G}$: 
the genus $g'$ of $Y$, the number $d$
of branch points, and the orders $m_1, \dots , m_d$
of the local monodromies. These invariants form the primary numerical type.

Once the primary numerical type is fixed, then the determination of the
  irreducible components of $\mathfrak M_{g;G}$ 
with a given primary numerical type
is, by Riemann's existence theorem, equivalent to the determination
of the orbits of the  group $Map (g',d) \times 
Aut ( G)$ on the set of possible monodromies
$\mu$. Here $Map (g',d)$ is the mapping class 
group of the curve $ Y \setminus \sB$, a curve of 
genus $g'$
with $d$ points removed.

Thus the general problem is to try to determine some finer numerical invariants
which determine these orbits (equivalently, the above irreducible components).

The secondary numerical type consists  of  the equivalence class
(for automorphisms of $G$) of the function which 
to each conjugacy class $\sC$ in
$G$ associates the number of branch points whose local monodromy lies
in the class $\sC$.

It was shown in \cite{cyclic} that the primary 
and secondary numerical type suffice to determine
the irreducible components $\mathfrak M_{g;G}$ in the case where $G$ is cyclic.

In this paper and its sequel we shall be 
concerned with the case where $G$ is a dihedral 
group $\Dn$.
In this case one can define the  numerical type, which is nothing else than the
primary and secondary numerical type
unless $n$ is even and the monodromy $\mu'$ onto 
the Abelianization $( \bZ / 2 \bZ)^2$ of $\Dn$
determines an unramified covering of $Y$.

We conjecture that each numerical type determines 
only one irreducible component,
and we present the proof here for the case 
$g'=0$; we have also a proof in the unramified 
case with $g' >0$.

Of course one can ask similar questions for more 
general groups, abelian groups should be 
relatively easy,
whereas more general solvable groups could lead to remarkable difficulties.

On the  opposite side, there is the case where 
$G$ is a simple group: for this case
we would like to call attention to the stability result of \cite{DT}.

The stability result of \cite{DT} states that, in the unramified
case (where   primary and secondary invariants boil down to  only  one invariant,  namely the 
genus $g'$), the number of irreducible components 
becomes
a constant independent of $g'$ for
$g'$ sufficiently large.

A very interesting question is whether a similar 
stability result holds fixing the secondary
numerical type but letting the genus $g'$ become sufficiently large.

\section{Preliminaries}

\paragraph{Dihedral groups.}
The \textit{dihedral group} $\Dn$ of order $2n$ 
is the group of symmetries of a regular
  polygon with $n$ edges.
We assume $ n \geq 3$, else we get the group $(\bZ/2 \bZ)^2$.

A simple representation of $\Dn$ is as the normal 
subgroup of  the affine group $A (1, \bZ / n\bZ )$
consisting of transformations of
the form
$$ m \mapsto \pm m + j ,  \ j \in \bZ / n\bZ .$$
It has generators $x$ such that $x (m) = m + 1$, 
and $y$ such that $y (m) = - m $;
$x^j$ corresponds to a rotation of $2\pi j /n$ 
around its barycenter and each element $x^j y$
( such that  $x^j y (m) = -m + j$) corresponds to 
a reflection with respect to a line of symmetry.

$\Dn$ can be defined by generators and relations as follows:
$$
\Dn=\langle\, x\, , \, y\, | \, x^n=y^2=(xy)^2=1\, \rangle \, .
$$
The above presentation shows right away that the 
Abelianization of $\Dn$ has the presentation
$$
\Dn^{Ab}=\langle\, x\, , \, y\, | \, xy = yx , 
y^2=(x)^{GCD(2,n)}=1\, \rangle \,
$$
hence we get $\bZ / 2\bZ $ for $n$ odd, $(\bZ / 2\bZ)^2  $ for $n$ even.

The $n$ \textit{reflections} $y\, , \, xy\, , \, 
\dots \, , \, x^{n-1}y$  will also be denoted
either by $s_0\, , \, s_1\, , \, \dots \, , \, 
s_{n-1}$ or by their indices $0\, , \, 1\, , \, 
\dots \, , \, n-1\, $.

For any  rotation $x^i$, its conjugacy class 
consists exactly of the elements $x^i$ and 
$x^{-i}$
(if $n=2i$ we obtain in this way the only central element).

If $n$ is odd all the
reflections  belong to the same conjugacy
class, while if $n$ is even two reflections 
$x^iy$ and $x^jy$ are conjugate if and only if
  $i \equiv j \, \mbox{(mod
$2$)}$.

These two cases are distinguished by the property 
of the corresponding affine transformation to
have fixed points, since  $x^iy (m) = m 
\Leftrightarrow  i \equiv 2 m (mod \ n) $, and 
this equation
has no solution if $n$ is even and $i$ is odd.

The \textit{automorphism group}  ${\rm Aut}(\Dn)$ 
is identified  with $A (1, \bZ / n\bZ )
\cong\bZ/n\bZ
\rtimes  (\bZ/n\bZ)^*$ as follows: the map $ 
\bZ/n\bZ\rightarrow {\rm Aut}(\Dn)$, which assigns
$i\in \bZ/n\bZ$ to the automorphism defined by 
$y\mapsto x^iy$ and $x\mapsto x$, identifies
$ \bZ/n\bZ$ with the normal subgroup of ${\rm 
Aut}(\Dn)$, consisting of those automorphisms
  which act trivially
on the subgroup of rotations.

  The quotient ${\rm Aut}(\Dn)/( \bZ/n\bZ)$ is
isomorphic to the  subgroup of  ${\rm Aut}(\Dn)$ 
consisting of automorphisms of the form $y\mapsto 
y$,
$x\mapsto x^i$ for
$i\in  (\bZ/n\bZ)^*$.

Observe that, $\Dn$ being a  normal subgroup of 
$A (1, \bZ / n\bZ )$, we get by conjugation
a homomorphism  $A (1, \bZ / n\bZ ) \ra {\rm 
Aut}(\Dn)$ which is an isomorphism exactly for 
$n$ odd.

\paragraph{Dihedral coverings of curves.} Let 
$Y$ be a compact connected Riemann surface of 
genus $g'$.
A \textit{dihedral covering} of $Y$ is a Galois 
covering $\pi \colon X\rightarrow Y$ with
  Galois group  $G= \Dn$
  and with $X$ connected. We will also  say that $\pi$ is a $\Dn$-covering.
\medskip\\
\textit{Riemann's existence theorem} allows us to 
use combinatorial methods to study 
$\Dn$-coverings, or
more generally any  $G$-covering (a Galois 
covering with an arbitrary finite Galois group 
$G$).

Let $\sB=\{\, y_1\, , \, \dots \, , \, y_d\, \}\subset Y$
be the branch locus of $\pi$.

Fix a base point $y_0\in Y\setminus \sB$ and a point $x_0\in \pi^{-1}(y_0)$.

Monodromy  gives a surjective group-homomorphism
\begin{equation}\label{m}
\mu \colon \pi_1 (\, Y\setminus \sB \, , \, y_0 \, ) \longrightarrow G \, .
\end{equation}

We recall now the definition of a {\bf geometric 
basis of $\pi_1 (\, Y\setminus \sB \, , \, y_0 \, 
)$}.

Let $\a_1\, , \, \b_1\, , \, \dots \, , \, 
\a_{g'}\, , \, \b_{g'}$ be simple nonintersecting 
(except in $y_0$)
closed arcs in $Y\setminus \sB$
which are based on $y_0$ and whose classes in 
$H_1(Y;\bZ)$ form a symplectic basis.

Let $\tilde{\cg}_i$ be an arc connecting $y_0$ with $y_i$,  contained in
$(Y\setminus \{\, \a_1\, , \, \b_1\, , \, \dots 
\, , \, \a_{g'}\, , \, \b_{g'}\, \})\cup \{ y_0 
\}$ and
such that $\tilde{\cg}_i$
intersects $\tilde{\cg}_j$ only in $y_0$ for $i\not=j$. Require moreover that
$\tilde{\cg}_1, \dots \tilde{\cg}_d$
stem out of $y_0$ with distinct tangents and 
following each other in counterclockwise order.

  Let $\cg_1\, , \, \dots \, , \, \cg_d
\subset Y\setminus (\sB \cup\{\, \a_1\, , \, 
\b_1\, , \, \dots \, , \, \a_{g'}\, , \, 
\b_{g'}\, \})
\cup \{ \, y_0\, \}$  be arcs defined as follows: 
$\cg_i$ begins at $y_0$, travels along 
$\tilde{\cg}_i$
to a point near $y_i$,  makes a small simple counterclockwise loop around $y_i$
and then returns to $y_0$ along
$\tilde{\cg}_i$.

Then we have chosen a geometric basis, and we have a presentation:
\begin{equation*}
\pi_1 (\, Y\setminus \sB\, , \, y_0 \, ) = 
\langle \, \a_1\, , \, \b_1\, , \, \dots \, , \, 
\a_{g'}\, ,
\, \b_{g'}\, , \, \cg_1\,
,
\, \dots \, , \, \cg_d\, | \, \Pi_{i=1}^{g'}[\, 
\a_i\, , \, \b_i\, ]\cdot \cg_1\cdot \, \dots \,
\cdot \cg_d=1\, \rangle\, .
\end{equation*}
Let $ \mathbf{T}(g', d)$ be the group defined 
abstractly by generators and relations as follows:
$$
\mathbf{T}(g', d):= \langle \, A_1\, , \, B_1\, , 
\, \dots \, , \, A_{g'}\, , \, B_{g'}\, , \,
  \Gamma_1\, , \, \dots \, , \, \Gamma_d\, | \, 
\Pi_{i=1}^{g'}[\, A_i\, , \, B_i\, ]\cdot
\Gamma_1\cdot \, \dots \,
\cdot \Gamma_d=1\, \rangle \, .
$$
The choice of a geometric basis yields
an obvious isomorphism
  $\pi_1 (\, Y\setminus \sB \, , \, y_0 \, )\cong \mathbf{T}(g', d)$ and
under this identification the homomorphism 
\eqref{m} corresponds to an epimorphism:
\begin{equation}\label{bfm}
\mbox{\boldmath${\mu}$} \colon  \mathbf{T}(g', d) \longrightarrow  G\, .
\end{equation}
Conversely, given a surjective homomorphism 
$\mbox{\boldmath${\mu}$} $ as in \eqref{bfm} such 
that
$\mbox{\boldmath${\mu}$} (\Gamma_i)\not= 1\ \forall i$,
  by Riemann's existence theorem the choice of a 
geometric basis as above ensures  the existence 
of a
$G$-covering $\pi\colon X\ra Y$ branched on $\sB 
$ and whose monodromy is $\mbox{\boldmath${\mu}$} 
$.

Varying a covering in a flat family with 
connected base, there are some numerical 
invariants which remain
unchanged, the first ones being the respective 
genera $g, g'$ of the curves $X$, $Y$, which are 
related by the
Hurwitz formula:

$$2 (g-1) = | G | [ 2 (g'-1) + \sum_i ( 1 - 
\frac{1}{m_i})], \ \ m_i : = ord (\mu (\cg_i) ). 
$$

Observe moreover that a different choice of the 
geometric basis changes the generators
$\cg_i$, but does not change their conjugacy classes (up to permutation),
hence another numerical invariant is provided by the number of elements
$\mu (\cg_i) $ which belong to a fixed conjugacy class in the group $G$.

We formalize these invariants through  the following definitions.

\begin{definition}\label{Hurwitz vector}
A \textbf{$G$-Hurwitz vector} is an ordered sequence
\begin{equation}\label{bfv}
\mathbf{v}=\left( \, a_1\, , \, b_1\, , \, \dots \, , \, a_{g'}\, , \, b_{g'}\,
, \, c_1\, , \, \dots \, , \,  c_d\, \right)
\in G^{2g' + d}
\end{equation}
such that the following conditions are satisfied:
\begin{description}
\item[(i)] $c_i\not=1$ for all $i$;
\item[(ii)] $G$ is generated by the components of 
$\vv$, $G=\langle \vv\, \rangle$;
\item[(iii)] $\Pi_{i=1}^{g'}[\, a_i\, , \, b_i\, 
]\cdot c_1\cdot \, \dots \, \cdot c_d=1$.
\end{description}
\end{definition}

To any $\Dn$-Hurwitz vector $\mathbf{v}$ we 
associate a tuple of positive integers $\nu 
(\mathbf{v})$
defined as follows.

If $n=2n' +1$ is odd, $\nu (\mathbf{v})=(k,k_1,\dots , k_{n'})$,
where $k$ (resp. $k_i$) is the number of the 
$c_i$'s in the conjugacy class of $y$ (resp. 
$x^i$).

If $n=2n'$ is even, $\nu 
(\mathbf{v})=(k_y,k_{xy},k_1,\dots , k_{n'})$, 
where $k_y$ (resp. $k_{xy}$, $k_i$)
  is the number of the $c_i$'s in the conjugacy 
class of $y$ (resp. $xy$, $x^i$).

The group  ${\rm Aut}(\Dn)$ acts diagonally on 
the set of Hurwitz vectors. This action induces 
an action of
${\rm Aut}(\Dn)$ on the set
$\mathcal{N}:=\{ \nu (\mathbf{v})\, | \, 
\mathbf{v}\, \mbox{is a Hurwitz vector}\, \}$ 
such that the map
$\nu$ is ${\rm Aut}(\Dn)$-equivariant.

The equivalence class of $\nu(\mathbf{v})$ in $\mathcal{N}/{\rm Aut}(\Dn)$
will be denoted by $[\nu(\mathbf{v})]$.

\begin{definition}\label{numerical type}
The \textbf{numerical type} of the Hurwitz vector $\vv$ is defined as follows.

If $n=2n' +1$ is odd, it is the pair $(g', 
[\nu(\mathbf{v})])$, where $g'$ is the genus of 
$Y$ and
$[\nu(\mathbf{v})]\in \mathcal{N}/{\rm Aut}(\Dn)$ is as above. \\
  If $n=2n'$ is even, then there are two cases:
\begin{itemize}
\item
Let $\epsilon \colon  D_{2n'} \ra (\bZ / 2 
\bZ)^2$ be the canonical surjection onto the 
Abelianization,
and let $p\colon Z \ra Y$ be the   degree 4 
covering associated to the composition
$\mu' : = \epsilon \circ \mu$.

Observe that $p$ is unramified if and only if 
none of the elements $c_i$ is a reflection
or a rotation $x^i$ with odd exponent $i$,

  By the Hurwitz
formula applied to $\mu$ and $\mu '$, the 
geometrical property that $p$ is unramified
is just a property of
$\mu$.
\item
(I) If $p$ is not unramified, then the numerical 
type  is again the pair $(g', [\nu(\mathbf{v})])$.
\item
(II) If $p$ is  unramified, then consider the two 
dimensional subspace $U$ of $H^1(Y, \bZ / 2 \bZ)$
dual to the surjection $H_1(Y, \bZ / 2 \bZ) \ra 
(\bZ / 2 \bZ)^2$ through which $\mu '$ factors.

Define $\iota \in \{0,1\}$ to be $=0$ if $U$ is isotropic, and $=1$ otherwise.

Then  the numerical type  is defined as the triple $(g', [\nu(\mathbf{v})],
\iota)$.\footnote{for $g'=0$ only case (I)
occurs.}
\end{itemize}

\end{definition}

\paragraph{Topological type.}

We recall a result contained in \cite{FabIso}, see also \cite{cime}.

Define the orbifold fundamental group 
$\pi_1^{orb} (\, Y\setminus \sB\, , \, y_0 \, ; 
m_1, \dots m_d ) $ of the
covering as

\begin{equation*}
  \langle \, \a_1\, , \, \b_1\, , \, \dots \, , \, \a_{g'}\, ,
\, \b_{g'}\, , \, \cg_1\,
,
\, \dots \, , \, \cg_d\, | \, \Pi_{i=1}^{g'}[\, 
\a_i\, , \, \b_i\, ]\cdot \cg_1\cdot \, \dots \,
\cdot \cg_d=1\, , \cg_j^{m_j} = 1 \ \forall j=1, \dots d\rangle\, .
\end{equation*}

We have then an exact sequence
  $$1 \ra \pi_1 (\, X \, , \, x_0 \,) \ra \pi_1^{orb}
(\, Y\setminus \sB\, , \, y_0 \, ; m_1, \dots m_d) \ra G \ra 1$$
which is completely determined by the monodromy, and which in turn determines,
via conjugation, a homomorphism

$$\rho \colon G \ra Out (\pi_1 (\, X \, , \, x_0 
\,)) = Map (X) : = Diff^+ (X) / Diff^0(X)  $$
which is fully equivalent to the topological action of $G$ on $X$.

We have that, by proposition 4.13 of 
\cite{FabIso}, all the curves $X$ of a fixed 
genus $g$
which admit a given topological action $\rho$ of 
the group $G$, specified up to an automorphism of 
$G$,
are parametrized by a connected complex manifold; 
arguing as in Theorem 2.4 of \cite{cyclic} we get

\begin{thm}

The triples $(X,G, \rho)$ where $X$ is a
complex projective curve of genus $g \geq
2$, and $G$ is a finite  group acting on $X$ with a topological action of
type $\rho$
  are  parametrized by a
connected complex manifold  $\sT_{g;G,\rho}$ of
dimension $3 (g'-1) + d $,  where $g'$ is the genus of $Y: = X / G$,
and  $d$ is the cardinality of the branch locus $\sB$.

The image  $\mathfrak M_{g;G,\rho}$ of
  $\sT_{g;G,\rho}$ inside the
moduli space
    $\mathfrak M_g$ is an irreducible closed 
subset of the same dimension  $3 (g'-1) + d $.

\end{thm}

The next question which the above result 
motivates is: when do two Galois monodromies
$\mu_1, \mu_2 :  \pi_1^{orb}
(\, Y\setminus \sB\, , \, y_0 \, ; m_1, \dots 
m_d) \ra G$ have the same topological type?

The answer is theoretically easy, since if the two covering spaces
have the same topological type then they are homeomorphic,
hence this means that the two monodromies differ by:

\begin{itemize}
\item
An automorphism of $G$.
\item
And a different choice of a geometric basis. This 
is performed by the mapping class group
(the first equality follows since the points of 
$\sB$ are the ends of $Y \setminus \sB$):
$$Map (Y, \sB) \cong  Map(Y \setminus \sB) := 
Diff^+(Y \setminus \sB) / Diff^0(Y \setminus \sB) 
.$$

\end{itemize}

\paragraph{Moduli spaces.}
Fixing a genus $g$ and a finite group $G$ we have a finite number of closed
irreducible subsets $\mathfrak M_{g;G,\rho} \subset \mathfrak M_g $
corresponding to the choice of a topological type $\rho$ for the action of $G$.

A first invariant for the topological type $\rho$ 
is provided by the pair $(g', d)$ where $g'$ is
as above the genus of $ Y : = X/G$, and $d$ is 
the cardinality of the branch locus $\sB \subset 
Y$.

A further numerical invariant is the $Aut(G)$ 
equivalence class of the class function
$\nu$ which, for  each conjugacy class
$\sC$ in
$G$, counts the number of local monodromies $c_i 
: = \mu(\cg_i)$ which belong to the conjugacy 
class $\sC$.

In particular, a weaker numerical invariant is 
given by the sequence of multiplicities $m_i$ of 
the branch points
($m_i$ is the order of $\mu (\cg_i$); one can 
assume w.l.o.g. $ m_1 \leq m_2 \leq \dots \leq 
m_d$).

One can consider then  the set of equivalence 
classes of  pairs $ (X,a)$, where $X$ is a 
projective curve of genus
$g$ and
$a$ is an effective action of $G$ on $X$ with primary numerical invariants
$(g; m_1, \dots m_d)$.

Two such pairs $(X,a)$ and $(X', a')$ are 
considered equivalent iff  there exists a 
biholomorphic map
$F\colon X\ra X'$ and an automorphism $\varphi\in {\rm Aut}(G)$ such that
  $F(gx)=\varphi(x)F(x)$, for any $x\in X$ and $g\in G$.

The set of such irreducible components $\mathfrak 
M_{g;G,\rho}$ with the given primary numerical 
invariants
$(g; m_1, \dots m_d)$ is then computed as the 
number of orbits of $Map (g',d) \times Aut (G) $
on the set of surjective homomorphisms
$$\mu :
\mathbf{T}(g', d; m_1, \dots m_d)
   \ra G  $$
where
$$\mathbf{T}(g', d; m_1, \dots m_d):=$$ $$:= 
\langle \, A_1\, , \, B_1\, , \, \dots \, ,
\, A_{g'}\, , \, B_{g'}\, , \,
  \Gamma_1\, , \, \dots \, , \, \Gamma_d\, | \, 
\Pi_{i=1}^{g'}[\, A_i\, , \, B_i\, ]\cdot
\Gamma_1\cdot \, \dots \,
\cdot \Gamma_d=1\,, \Gamma_i^{m_i} =1 \ \forall i \rangle \, .$$

The geometrical insight is that the union of such 
components $\mathfrak M_{g;G,\rho}$ has a finite 
map $Q \colon \mathfrak M_{g;G,\rho} \ra \mathfrak M_{g', d}$
onto the (coarse) moduli space  $\mathfrak M_{g', d}$ of smooth curves
of genus $g'$ with $d$ unordered marked points. 
This is a topological covering and
the fundamental group of the base is a quotient 
of the mapping class group  $Map (g',d)$.

Hence the components $\mathfrak M_{g;G,\rho}$ are 
detected by the orbits of the monodromy
of this covering space.

\paragraph{The case of the dihedral group.}

Let $n$ be a positive integer $n\geq 3$ and let
$(g, [\nu(\mathbf{v})])$ (resp. $(g, 
[\nu(\mathbf{v})], \iota)$ be a numerical type.

Let $\hfami_{\Dn}(g, [\nu(\mathbf{v})]) $ (resp. 
$\hfami_{\Dn}(g, [\nu(\mathbf{v})], \iota) $)
be the set of equivalence classes of  pairs $ 
(X,a)$, where $X$ is a Riemann surface of genus 
$g$ and
$a$ is an effective action of $\Dn$ on $X$ such 
that  the $\Dn$-covering $X\ra X/\Dn$ is of 
numerical type
$(g, [\nu(\mathbf{v})])$ (resp. $(g, [\nu(\mathbf{v})], \iota)$.\\

The main question we address is whether the 
spaces $\hfami_{\Dn}(g, [\nu(\mathbf{v})]) $, 
respectively
$\hfami_{\Dn}(g, [\nu(\mathbf{v})], \iota) $ are irreducible, i.e., are spaces
$\mathfrak M_{g;\Dn,\rho}$ for a unique topological type $\rho$.
This can be proved by showing the transitivity of $Map (g',d) \times Aut (\Dn)$
on the set of monodromies of given (full) numerical type.

This is the same thing as bringing each monodromy 
with a given numerical type to a normal form.

\section{The   case $g'=0$}

In this Section we assume $g'=0$.

The moduli space  $\mathfrak{M}_{0,d}$ is a quotient 
of  $( \frak{S}^d\bP^1)\setminus \Delta$ by the action of the projective linear group $\bP GL ( 2 , \bC )$,
where $\frak{S}^d\bP^1$ is  the $d$-th symmetric product  of $\bP^1$,
and $\Delta $ is the subset of $\frak{S}^d\bP^1$ 
consisting of points with two or more equal 
coordinates.

We have  $\frak{S}^d\bP^1\cong \bP^d$, therefore we consider  the action
of the \textit{braid group of the sphere}  $\sS \mathcal{B}_d:=
\pi_1 (\, \bP^d\setminus \Delta\, , \, \ul{y}\, )$
on the fibre  over $\ul{y}$ of the above map $Q$.

The group $\sS \mathcal{B}_d$ is a quotient of 
Artin's braid group $\mathcal{B}_d$
which is  generated by the
so-called
\textit{elementary braids} $\s_1\, , \, \dots \, 
, \, \s_{d-1}$ acting on the Hurwitz vector $\vv=(\, 
c_1\, , \, \dots \, , \, c_d\,)$ as follows:
\begin{eqnarray*}
(\, c_1\, , \, \dots \, , \, c_d\, )\s_i &=& (\, 
c_1\, , \, \dots \, , \, c_i c_{i+1}c_i^{-1}\, , 
\,  c_i\, , \, \dots \, , \, c_d\, ) \, , \\
(\, c_1\, , \, \dots \, , \, c_d\, )\s_i^{-1} &=& 
(\, c_1\, , \, \dots \, , \, c_{i+1}\, , \, 
c_{i+1}^{-1}c_i c_{i+1}\, , \, \dots \, , \, 
c_d\, )\, .
\end{eqnarray*}
Recall moreover the diagonal action of the  group ${\rm Aut}(\Dn)$   on the set of Hurwitz vectors.

Since the two actions commute, we have an action 
of the group $\mathcal{B}_d\times {\rm Aut}(\Dn)$.
\begin{definition}
Two Hurwitz vectors $\mathbf{v}$ and $\mathbf{w}$ are
said to be Hurwitz equivalent, or 
Braid-equivalent (resp. automorphism-equivalent,
braid-automorphism-equivalent)  if there exist 
$\s\in \mathcal{B}_d$ (resp. $\varphi \in {\rm 
Aut}(\Dn)$,
$(\s, \varphi)\in \mathcal{B}_d\times {\rm 
Aut}(\Dn)$) such that $\mathbf{w}=\mathbf{v}\s$ 
(resp.
$\mathbf{w}=\mathbf{v}\varphi$, 
$\mathbf{w}=\mathbf{v}(\s, \varphi)$). In this 
case we write
$\mathbf{v}\beq \mathbf{w}$ (resp. 
$\mathbf{v}\aeq \mathbf{w}$, $\mathbf{v}\baeq 
\mathbf{w}$).
\end{definition}

\begin{notation}\normalfont
Identify a reflection $s_i (m) = - m + i$ with 
its index $i \in \bZ / n \bZ$.
\end{notation}

The main result of this section is the following
\begin{thm}\label{genus 0}
The group $\mathcal{B}_d\times {\rm Aut}(\Dn)$ 
acts transitively on the set of Hurwitz vectors of
a fixed numerical type, hence dihedral covers of 
$\bP^1$ of a fixed numerical type
form an irreducible closed subvariety of the moduli space.

  More precisely, given $\vv$ with 
$\nu(\mathbf{v})=(k,k_1,\dots , k_{n'})$ (resp.
$\nu(\mathbf{v})=(k_y,k_{xy},k_1,\dots , k_{n'})$),
set $R: =\sum_i k_i$, and assume (w.l.o.g.)  $\{ 
h,k \}=\{ k_y,k_{xy} \}$, $h\leq k$
(observe that $k$, resp. $k+h$ is even).

We have then, assuming  throughout
$0< r_i \leq r_{i+1} \leq n'$,  $\ul{r}=(r_1 , \dots , r_R)$ and
setting

$|\ul{r}|: \equiv \sum_i r_i (mod \ n)$:
\begin{description}
\item[i)]
$\vv \baeq (\, \underbrace{0\, , \, \dots \, , \, 0\, , \, 1\, ,
\, 1+|\ul{r}| }_{k}\, , \, x^{r_1}\, , \, \dots 
\, , \, x^{r_R}\, )$,  if $n=2n'+1$.

\item[ii)]
$\vv \baeq (\, \underbrace{0\, , \, \dots \, , \, 
0}_{h}\, , \, \underbrace{1\, ,
\,\dots \, ,1\, , \, \lambda }_{k}\, , \, 
x^{r_1}\, , \, \dots \, , \, x^{r_R}\, )$,
  if $n=2n'$ and $h\not=0$.

Here $\lambda = |\ul{r}|  + \epsilon$, where 
$\epsilon \in \{0,1\} $,  $\epsilon + k \equiv 
1(mod \ 2).$
\item[iii)]
$\vv \baeq (\,  \underbrace{1\, , \,\dots \, ,1\, , \, 3\, , \,  \lambda }_{k}\,
  , \, x^{r_1}\, , \, \dots \, , \, x^{r_R}\, )$,  if $n=2n'$ and $h =0$.

Here $\lambda = |\ul{r}|+3 $.\end{description}

\end{thm}

  We collect in the next section some preliminary 
results that shall be used in the proof.

\begin{remark}
It was brought to our attention after the paper was completed that a rather complicated but more general classification of Hurwitz orbits 
on $\Dn^d$ was done in \cite{Sia}. It is however not clear to us  whether  one can deduce our theorem above from these results.
\end{remark}

\subsection{Auxiliary results}

\begin{remark}\normalfont
Identifying a reflection $s_i (m) = - m + i$ with 
its index $i \in \bZ / n \bZ$,
then  conjugation corresponds to the  action of another reflection on $i$:
$$
s_i \mapsto s_j s_i s_j \quad\text{corresponds to}\quad
i \mapsto 2j-i = j - (i-j).
$$
\end{remark}
\begin{remark}\normalfont
The action of  $\s_1$ on a pair of reflections 
$(i,j)$ leaves their product invariant, hence
leaves the difference $i-j$ invariant: for instance $(i,j)\s_1 = (2 i - j, i)$.

\end{remark}

\begin{lemma}[Normalization of reflection triples]
\label{triple}
Given a sequence of reflections $(i,j,k)$ in $D_n$ which generate
a dihedral subgroup $D_m$, its Hurwitz orbit contains
a sequence of type $(i',j',j')$ and a sequence of type $(i'',i'',j'')$.

In particular the subgroup $D_m$ is generated by the first two
entries of a suitable sequence in the Hurwitz orbit.
\end{lemma}

\proof
First we consider the action of the four elements $\s_1,\s_2,\s_1\inv,
\s_2\inv$ on the triple.
$$
\begin{array}{c}
(i,j,k)\s_1 = (2i-j, i,k),\quad
(i,j,k)\s_2 = (i,2j-k,j),\\
(i,j,k)\s_1\inv = (j,2j-i,k),\quad
(i,j,k)\s_2\inv = (i,k,2k-j).
\end{array}
$$
The corresponding transformations on the 
differences of consecutive elements are
$$
\begin{array}{c}
(j-i,k-j)\s_1 = (j-i,(k-j)+(j-i)),\quad
(j-i,k-j)\s_2 = ((j-i)-(k-j),k-j),\\
(j-i,k-j)\s_1\inv = (j-i,(k-j)-(j-i)),\quad
(j-i,k-j)\s_2\inv = ((j-i)+(k-j),k-j).
\end{array}
$$
As long as both differences are non-zero (modulo $n$), we can reduce
the maximal difference by one of these transformations.

This process must terminate, hence we reach a situation
where one of the differences is zero.

We can arrange for the other difference to become zero
by at most two additional transformations.
Then we end up with a triple such that the last two
entries are equal, and also with a triple such that the first two
entries are equal.

Note that we can compute the necessary transformations
using the Euclidean algorithm for  the two differences.
\qed

\begin{lemma}\label{abm}
Let $(s_i,s_j,x^m)\in \Dn^3$. Then,  for all $\ell \in \bZ$ we have:

$(\, s_i \, , \, s_j \, , \, x^m\, )\beq (\, 
s_{i+2\ell m}\, , \, s_{j+2\ell m} \, , \, x^m\, 
)$.
\end{lemma}
\proof
For any $\ell \in \bZ$ the following formula can be verified:
\begin{equation}\label{abm1}
(\, s_i\, , \, s_j\, , \, x^m\, )(\s_2 \s_1 \s_1 
\s_2)^{\ell}=(\, s_{i+2\ell m}\,
, \, s_{j+2\ell m}\, , \, x^m\, )\, \nonumber .
\end{equation}
This proves the claim.
\qed

\begin{lemma}[Double exchange]\label{triple}
\label{repetition}
The following equivalence holds:

$ (j, i, i) \quad \beq \quad (i,i,j). $

\end{lemma}
\proof
This follows from the following equality:
$$
(j,i,i)\sigma_1^{-1}\sigma_2^{-1}=(i,2i-j ,i)\sigma_2^{-1} =(i,i,j).
$$
\qed

\begin{lemma}[Normalization of pair sequences]
\label{repetition}
The following equivalences hold.
\begin{enumerate}
\item
$ (0,0, i, i) \quad \beq \quad (0,0, -i, -i), $
\item
$ ( i,i,j,j) \quad \beq \quad (j,j, i,i), $
\item
$ ( i, i, j, j) \quad \beq \quad (\, i+\ell (j-i) 
\, ,  \, i+\ell (j-i) \, , \, j+\ell (j-i) \, ,
  \, j+\ell (j-i)\, ) $, $\forall \ell \in \bZ$.
\item
$ (0,0, i,i,j,j)\quad \beq \quad(0,0,i,i,j-2\ell 
i,j-2\ell i) $, for any $\ell\in \bN$
\end{enumerate}
\end{lemma}

\proof
We achieve equivalence (i) by
\begin{eqnarray*}
(0,0,i,i) \s_2\s_3^2\s_2 & = & (0, -i,-i,0) \s_3\s_2
\\
& = & (0,-i,-2i,-i)\s_2 \\
& = & (0,0,-i,-i).
\end{eqnarray*}
We achieve equivalence (ii) by applying twice  Lemma \ref{triple}.\\
\\
Equivalence (iii) follows from the formula
$$
(\, i\, , \, i \, , \, j\, , \, j\, )(\s_2^{-1} 
\s_1^{-1}\s_3 \s_2)^{\ell}= (\, i+\ell (j-i) \,
  ,  \, i+\ell (j-i) \, , \, j+\ell (j-i) \, , \, 
j+\ell (j-i)\, ) \, , \quad \forall \ell \in \bZ 
\, ,
$$
which can be proved e.g. by induction.\\

For  equivalence (iv) we  have:
\begin{eqnarray*}
(0,0,i,i,j,j)\s_4 \s_5 \s_5 \s_4&=& (0,0,i,2i-j , 2i-j , i) \s_5 \s_4\\
&=& (0,0,i,i,2i-j , 2i-j ) \\
&\beq & (0,0,i,i,j-2i , j-2i )\quad \mbox{by (i)}\, .
\end{eqnarray*}
Iterating this procedure we get the claim for all $\ell \in \bN$.

\qed

\begin{lemma}\label{no rotations}
Let $n$ be an integer $n\geq 3$ and let  $\vv=(\, 
i_1\, , \,  \dots \, , \,  i_{2N}\, ) \in 
\Dn^{2N} $
  be a Hurwitz vector  whose components are all reflections.
Then there exists   $j\in \bZ / n\bZ$ such that:
$$
\vv \baeq \begin{cases} (\, 0 \, , \, \dots \, , 
\,  0 \, , \,  j \, , \,  j \, )  & \mbox{if $n$ 
is odd};\\
(\, 0 \, , \, \dots \, , \,  0 \, , \,  j \, , \, 
\dots \, , \, j \, ) & \mbox{if $n$ is even} \, . 
\end{cases}
$$
Moreover, the automorphisms involved in the 
previous equivalences are all of the form
$y\mapsto x^\ell y$, $x\mapsto x$.
\end{lemma}

\proof
Using Lemma \ref{triple} inductively we get
$$
\vv \beq (i_1,i_1,i_2,i_2,\dots,i_{N-1},i_{N-1},i_N,j_N).
$$
Then also $i_N=j_N$, since the product is the 
identity and we have in fact obtained
$$
\vv  \beq (i_1,i_1,i_2,i_2,\dots,i_{N-1},i_{N-1},i_N,i_N).
$$
By the automorphism $y\mapsto x^{-i_1} y$, $x\mapsto x$ which is of the form
given in the claim of the lemma we get the following form
$$
\baeq (0,0, i_2,i_2, \dots, i_N,i_N)
$$
and we may assume $i_\nu\geq0$ by choosing suitable representatives.

The assertion of the Lemma follows now  from the following:

CLAIM : Unless the sequence is already in the  asserted form
there is another sequence of pairs of non-negative integers representing
a Hurwitz vector in the same braid equivalence class
which is strictly smaller with respect to the lexicographical ordering.

Since we may reorder the sequence of pairs according to
Lemma \ref{repetition} $ii)$, we may assume
$0\leq i_2\leq \dots \leq i_N$.

Assume now we have three different kinds of 
entries $0<i<j$. Then by using once more
Lemma \ref{repetition} $ii)$ we can bring these entries next to each other
and have then a subsequence of the form $(0,0,i,i,j,j)$.

By Lemma \ref{repetition} $iv)$ and $i,ii)$ we have:
\begin{eqnarray}
(0,0,i,i,j,j) & \beq & (0,0,i,i,j-2i,j-2i) \\
& \beq & (0,0,i,i, -j+2i, -j+2i)
\end{eqnarray}
Now, either $j-2i\geq 0$ or $j-2i<0$.
In the first case $j>j-2i\geq0$ and the r.h.s.\ 
of $(4)$ is smaller than the l.h.s.
In the second case $0<i<j$ implies $2i-j<j$ and the r.h.s.\ of $(5)$
is smaller than the l.h.s\ of $(4)$.

Therefore we can reduce to a sequence of pairs where all entries are
either $0$ or a positive integer $j$.
This concludes the claim in case where $n$ is even.

In the case where $n$ is odd we may have reached a situation with at least
four entries equal to $j$. But, according to Lemma \ref{repetition} $iv)$
with $\ell=-(n-1)/2$, we have
\begin{eqnarray*}
(0,0,j,j,j,j) & \beq & (0,0,j,j,j-2\ell j,j-2\ell j) \\
& = & (0,0,j,j,nj,nj) \\
& = & (0,0,j,j,0,0) \\
& \beq & (0,0,0,0,j,j)
\end{eqnarray*}
Hence also in this case our claim holds true.
\qed

\begin{lemma}[Normalization of reflection pair]
\label{pair}
Given a sequence of reflections $(i_0,j_0)$ in $D_n$ which generate
a dihedral subgroup $D_m$, its Hurwitz orbit consists of
the pairs $(i,j)$ with $i\equiv i_0 \left( \mod \frac n m \right)$ and $j-i=j_0-i_0$.
\end{lemma}

\subsection{Proof of Theorem \ref{genus 0}}
\noindent 1. We first bring all the rotations to 
the right by  elementary braids,  obtaining:
\begin{equation}\label{1}
\vv\beq (\, s_{i_1}\, , \, \dots \, , \, 
s_{i_{2N}}\, , \, x^{r_1}\, , \, \dots \, , \, 
x^{r_R}\, )\, ,
\end{equation}
where $2N=k$ if $n$ is odd, $2N=h+k$ if $n$ is even and $R=\sum_i k_i$.

Observe moreover that we can arbitrarily permute 
the rotations among themselves,
a fact that at a later moment will allow us to 
assume $ r_i \leq r_{i+1}$, $\forall i$.
\medskip\\
2. If $r_i > n'$,  we bring $r_i$ next to the 
reflection $s_j = s_{i_{2N}}$ and then
apply  a full twist with this reflection $s_j$, obtaining:
$$
(\, s_j\, , \, x^{r_i}\, )\s_1^2 = (\, s_{j-2r_i}\, , \, x^{-r_i}\, )\, .
$$
Hence we can assume $0<r_i\leq n'$ for all $i$.
\medskip\\
3. If $n$ is even, without loss of generality, we may further assume that
$$k=| \{\, s_i\, | \, s_i\, \mbox{is conjugate to} \, s_{i_{2N}}\, \} | \, .$$
\medskip\\
4. By Lemma \ref{triple} we can assume that $i_{2f} = i_{2f-1}$ for any
  $f\in \{ \, 1,\dots , N-1\, \}$.

Then we set $j_f=i_{2f}$ for
$f\in \{ \, 1,\dots , N-1\, \}$ and $j_N= i_{2N-1}$, thus we have:
$$
\vv \beq (\, j_1\, , \, j_1\, , \, j_2\, , \, j_2\, , \, \dots \, , j_{N-1}\,
  , \, j_{N-1} \, , \, j_N\, , \, j_N + |\ul{r}| 
\, , \, x^{r_1}\, , \, \dots \, , \, x^{r_R} \, 
)\, .
$$
Notice that   the condition that $k=| \{\, s_i\, 
| \, s_i\, \mbox{is conjugate to}
  \, s_{j_N+ |\ul{r}|}\, \} | $ still
holds.
\medskip\\
5. Consider the vector
$$
\mathbf{w}:=(\, j_1\, , \, j_1\, , \, j_2\, , \, 
j_2\, , \, \dots \, , j_N\, , \, j_N\, ) \in 
\Dn^{2N}\, .
$$
The subgroup $\langle \mathbf{w}\rangle\leq \Dn$ generated by $\mathbf{w}$ is
isomorphic either to $\bZ/2\bZ$,  $\bZ/2\bZ \times \bZ/2\bZ$
or to $\Dm$ with $m\geq 3$.

We  show now, in each of these three cases, that 
$\vv$ is equivalent to one of the vectors in the 
statement of
Theorem
\ref{genus 0}.
\medskip

\noindent \ul{I. $\langle \mathbf{w}\rangle \cong 
\bZ/2\bZ$}. Then  $j_1=j_2=\dots =j_N=:j$ and
$$
\mathbf{v} \beq (\, j\, , \, \dots \, , \, j \, , 
\, j+|\ul{r}|\, , \, x^{\ul{r}}\, )\, .
$$
We have that $\Dn=\langle \, j\, , \, 
x^{\ul{r}}\, \rangle$ and hence $GCD(n,\ul{r})=1$.
By  Lemmas \ref{abm} and \ref{repetition} ii) it follows that
$$
\mathbf{v}\beq (j+2\ell\, , \, \dots \, , \, 
j+2\ell \, , \, j+2m\, , \, j+2m +|\ul{r}|\,
  , \, x^{\ul{r}})\, , \, \forall \ell \, , \, m \in \bZ\, .
$$
If $n$ is odd, $2$ is a generator of $\bZ/n\bZ$ and hence the result follows.

If $n$ is even, we may assume that $j$ is odd. In the case where
moreover $|\ul{r}|$ is even (i.e., $h=0$)
we choose $\ell$ such that $j+2\ell =1$ and   we 
set $m=(3-j)/2$. Otherwise we take $m$ such that
$j+2m+|\ul{r}|=0$, and $\ell$ such that $j+2\ell 
=- 1$ . The result follows by a sequence of 
Hurwitz moves
between reflections bringing the element $0$ from 
the last position to the initial one.

\medskip

\noindent \ul{II. $\langle \mathbf{w}\rangle\cong \bZ/2\bZ \times \bZ/2\bZ $}.

  Then  $n$ must be even and
there exist $i\, , \, j \in \Dn$ with $i - j = 
\frac{n}{2}$ such that $j_f\in \{\, i\, , \, j \, 
\}$
for all $f\in \{1, \dots , N\}$.
Using Lemma \ref{repetition} (ii) we can bring 
all the pairs of $i$'s to the left and obtain:
\begin{eqnarray}
\mathbf{v} &\beq& (\, i \, , \, \dots \, , \,i 
\, , \, j \, , \, \dots \,  , \, j \, , \, j 
+|\ul{r}|\,
  , \, x^{\ul{r}}\, )\quad \mbox{or} \label{z2xz2 1} \nonumber \\
\vv&\beq& (\, i \, , \, \dots \, , \,i  \, , \, j 
\, , \, \dots \,  , \, j \, , \, i\, , \,  i 
+|\ul{r}|\,
  , \, x^{\ul{r}}\, )\, . \nonumber  
  \label{z2xz2 2}\,
\end{eqnarray}

Exchanging the roles of $i$ and $j$ and using 
again   Lemma \ref{repetition} (ii) we see that 
the second vector
is braid-equivalent to one of the first type, hence
we shall only consider the first vector.

  We have that $\Dn = \langle \, i \, , \, j\, , 
\, x^{\ul{r}}\, \rangle $ and so 
$GCD(n,\ul{r})\in \{1,2\}$.

  If $GCD(n,\ul{r})=1$,   we apply Lemma \ref{abm} and we get:
  \begin{equation}\label{z2z21}
  \mathbf{v} \beq  (\, i +2\ell \, , \, \dots \, , 
\,i +2\ell  \, , \, j +2m \, , \, \dots \,  , \, 
j+2m \,
  , \,j+2p\, , \,  j+2p +|\ul{r}|\, , \, 
x^{\ul{r}}\, ) \, , \, \forall \ell, m, p\in 
\bZ\, . \nonumber
  \end{equation}
For an appropriate choice of $\ell, m, p\in \bZ$ 
if $i,j$ are odd we reach the required normal 
form iii).

If one is even and the other is odd there are two possibilities:
either the larger group of $k$ elements (to which $ j + |r|$ by our assumption belongs)
contains the numbers $j$, or that it contains the numbers $i$.

In the former case with an 
automorphism of $\Dn$ we achieve that $i$ is even
and again for  an appropriate choice of $\ell, m, 
p\in \bZ$  we reach the required normal form ii).

In the latter case since the number of occurrences of $i$ is even we first apply repeatedly Lemma \ref{triple} to move all the $j$'s to the left,
then with an automorphism we achieve that $j$ is even, and finally for an appropriate 
choice of $\ell , m$  we reach the required normal form ii).

  Assume now that $GCD(n,\ul{r})=2$.

Then $\frac{n}{2}$ must be odd
  and therefore $i$ and $j$ have different 
parities. Moreover $|\ul{r}|$ is even and so we
may assume, acting with an automorphism, that $h$
coincides with the number of $i$'s.
  By Lemma \ref{abm} we get:
\begin{equation}\label{z2z22}
\mathbf{v} \beq  (\, i +4\ell \, , \, \dots \, , 
\,i +4\ell  \, , \, j \, , \, \dots \,
  , \, j \, , \,  j +|\ul{r}|\, , \, x^{\ul{r}}\, 
) \, , \, \forall \ell\in \bZ\, . \nonumber
  \end{equation}
If $i\equiv j-1 \, \mbox{(mod $4$)}$, we apply the automorphism
$x^{j-1}y \mapsto y$, $x\mapsto x$ to transform the vector 
in the desired form ii).
Otherwise $i\equiv j+1 \, \mbox{(mod $4$)}$ and we proceed as follows:
\begin{eqnarray*}
\mathbf{v} &\beq& (\, j+1 \, , \, \dots \, , 
\,j+1  \, , \, j \, , \, \dots \,  , \, j \,
  , \,  j +|\ul{r}|\, , \, x^{\ul{r}}\, ) \\
&\aeq& (\, -j-1 \, , \, \dots \, , \, -j-1  \, , 
\, -j \, , \, \dots \,  , \, -j \, , \,  -j 
-|\ul{r}|\, , \, x^{-\ul{r}}\, )\\
&\aeq& (\, 0 \, , \, \dots \, , \,0  \, , \, 1 \, 
, \, \dots \,  , \, 1 \, , \,  1 -|\ul{r}|\, , \, 
x^{-\ul{r}}\, )\\
&\beq& (\, 0 \, , \, \dots \, , \,0  \, , \, 1 \, 
, \, \dots \,  , \, 1 \, , \,  1 +|\ul{r}|\, , \, 
x^{\ul{r}}\, )\, ,
\end{eqnarray*}
where in the second equivalence we have used the 
automorphism $x\mapsto x^{-1}$, $xy\mapsto 
x^{-1}y$,
in the third equivalence we have used the 
automorphism $x\mapsto x$, $xy\mapsto x^{j+2}y$,
and in the fourth one we proceeded as in the reduction step 2.

\medskip

\noindent \ul{III. $\langle \mathbf{w}\rangle 
\cong \Dm$, $m\geq 3$}. 

By Lemma \ref{no rotations}  applied to $ \mathbf{w}$ we reduce $ \mathbf{w}$ to the form 
\begin{eqnarray}
\mathbf{w} &\baeq& (\, 0 \, , \, \dots \, , \, 0 
\, , \, j   \, , \, \dots \, , \, \, j  \, , \, 
\, j \, ).
\end{eqnarray}

Since we want to apply the corresponding moves to $\mathbf{v}$ we avoid moves which put the last pair into a different position.

By a careful modification of the proof of Lemma \ref{no rotations} this restriction leads to 

\begin{eqnarray}
\mathbf{w} &\baeq& (\, 0 \, , \, \dots \, , \, 0 
\, , \, j   \, , \, \dots \, , \, \, j  \, , \, 
\, j \,,  j \,) \quad \mbox{or}  \\
&\baeq& (\, 0 \, , \, \dots \, , \, 0 \, , \, j 
\, , \, \dots \, , \, \, j  \, , \, \, 0, 
0 \, ) \nonumber \, 
\end{eqnarray}

 hence we have
\begin{eqnarray}
\mathbf{v} &\baeq& (\, 0 \, , \, \dots \, , \, 0 
\, , \, j   \, , \, \dots \, , \, \, j  \, , \, 
\, j + |\ul{r}|\,
  , \, x^{\ul{r}}\, )\label{dm 1}\quad \mbox{or} \label{dm} \\
&\baeq& (\, 0 \, , \, \dots \, , \, 0 \, , \, j 
\, , \, \dots \, , \, \, j  \, , \, \, 0, 
|\ul{r}|\,
  , \, x^{\ul{r}}\, )\label{dm }\nonumber \, .
\end{eqnarray}
It is clearly enough to consider only the first case.

Observe that  $\Dn=\langle \, y \,
, \, x^jy \, , \, x^{\ul{r}}\, \rangle$. \\
If $n$ is odd, then by Lemma \ref{no rotations} we have:
$$
\mathbf{v} \, \baeq \,  (\, 0\, , \, \dots \, , 
\, 0 \, , \, j \, , \, j+ |\ul{r}|\, , \, 
x^{\ul{r}}\, )
$$
with $GCD(j,n, \ul{r})=1$.

From Lemma \ref{abm} it follows that the right hand side is braid-equivalent to
$$
(\, 0\, , \, \dots \, , \, 0 \, , \, j +2\ell M 
\, , \, j+ |\ul{r}| +2\ell M \, , \, x^{\ul{r}}\, 
)\, ,
$$
where $M:=GCD(n,\ul{r})$. We have then  $GCD(j,M)=1$.

Set $\ell=n/\gamma$, with $\gamma$ equal to the product
  of all common prime factors of $n$ and $j$ taken 
with the maximal power with which they
divide $n$; hence $GCD(\epsilon:= j +2\ell 
M,n)=1$. In fact if $ p | n$, then $p\neq 2$, and 
either $p| \ell$, or $p
| j$: but  if $p| \ell$ then $p | j$, contradicting 
that $\ell$ and $j$ are relatively prime; if instead 
$p | j$,
then by the same token $p|M$, contradicting $GCD(j,M)=1$.

Using Lemma \ref{pair} we get:
\begin{eqnarray}\label{n odd}
(\, 0\, , \, \dots \, , \, 0 \, , \, \epsilon \, , \, \epsilon + |\ul{r}|  \,
  , \, x^{\ul{r}}\, )&\beq& (\, 0\, , \, \dots \, 
, \, 0 \, , \, \ell \epsilon \,
  , \, (\ell +1) \epsilon \, , \, \epsilon +
|\ul{r}|  \, , \, x^{\ul{r}}\, )\nonumber \\
\mbox{(for $\ell = - \epsilon^{-1}$)} & = & (\, 0\, , \, \dots \, , \, 0 \,
  , \, -1 \, , \,  \epsilon -1 \, , \, \epsilon + 
|\ul{r}|  \, , \, x^{\ul{r}}\, )\nonumber \\
&\beq& (\, 0\, , \, \dots \, , \, -\ell \, , \, 
-\ell -1 \, , \,  \epsilon -1 \,
, \, \epsilon + |\ul{r}|  \, , \, x^{\ul{r}}\, ) \nonumber \\
\mbox{(for $\ell =- \epsilon$)} &=&(\, 0\, , \, \dots \, , \, \epsilon \,
  , \, \epsilon -1 \, , \,  \epsilon -1 \, , \, 
\epsilon + |\ul{r}|  \, , \, x^{\ul{r}}\, ) 
\nonumber
  \, . \nonumber
\end{eqnarray}
Repeating these steps inductively we obtain a vector of the following form:
\begin{eqnarray*}
(\,  \epsilon \, , \, \epsilon -1 \, , \,   \dots 
\, , \, \epsilon -1 \, , \, \epsilon + |\ul{r}| 
\,
  , \, x^{\ul{r}}\, ) &\beq& (\,  \epsilon -1  \, 
, \,   \dots \, , \, \epsilon -1 \,
, \, \epsilon  \, , \, \epsilon + |\ul{r}|
\, , \, x^{\ul{r}}\, )\\ &\aeq&  (\,  0  \, , \,   \dots \, , \, 0 \, , \, 1 \,
, \, 1 + |\ul{r}|  \, , \, x^{\ul{r}}\, )\, .
\end{eqnarray*}
The case where $n$ is odd  is then settled.

Let now $n$ be even.

We distinguish three cases: $h=0$, $h$
  is equal to the number of $y$'s in \eqref{dm}, 
or $h$ is equal to the number of $x^jy$'s in 
\eqref{dm}.

In the first case we apply Lemma \ref{abm} to obtain:
$$
  (\, 0 \, , \, \dots \, , \, 0 \, , \, j   \, , 
\, \dots \, , \, \, j  \, , \, \, j + |\ul{r}|\,
  , \, x^{\ul{r}}\, )\beq  (\, 0 \, , \, \dots \, , \, 0 \, , \, j+2\ell M   \,
, \, \dots \, , \, \, j +2\ell M \, , \, \, j +2\ell
M+ |\ul{r}|\, , \, x^{\ul{r}}\, )\, ,
$$
where again $M=GCD(n,\ul{r})$. Let $\ell=n/\gamma$, with $\gamma$
  equal to the product of all common prime factors 
of $n$ and $j$ taken with the maximal power with 
which
they divide $n$; hence $GCD(\epsilon:= j +2\ell M,n)=2$.

Using Lemma \ref{pair} we have:
\begin{eqnarray*}
(\, 0\, , \, \dots \, , \, 0 \, , \, \epsilon \, , \, \dots \, , \, \epsilon \,
  , \,   \epsilon+ |\ul{r}|\, , \, x^{\ul{r}}\, ) &\beq&
  (\, 0\, , \, \dots \, , \, 0 \, , \, \ell 
\epsilon \, , \,  (\ell +1) \epsilon \,
, \, \epsilon \, , \, \dots \, , \, \epsilon \, , 
\,   \epsilon+ |\ul{r}|\, , \, x^{\ul{r}}\, ) \\
  \mbox{(for $\ell \epsilon = -2$)} &=&  (\, 0\, , 
\, \dots \, , \, 0 \, , \, -2 \,
  , \,   \epsilon -2 \, , \, \epsilon \, , \, \dots \, , \, \epsilon \,
  , \,   \epsilon+ |\ul{r}|\, , \, x^{\ul{r}}\, ) \\
&\beq& (\, 0\, , \, \dots \, , \, -2\ell \, , \, -2(\ell +1) \,
  , \,   \epsilon -2 \, , \, \epsilon \, , \, \dots \, , \, \epsilon \,
  , \,   \epsilon+ |\ul{r}|\, , \, x^{\ul{r}}\, ) \\
  \mbox{(for $-2\ell = \epsilon $)} &=& (\, 0\, , \, \dots \,
, \, \epsilon \, , \, \epsilon -2 \, , \,   \epsilon -2 \, , \, \epsilon \,
  , \, \dots \, , \, \epsilon \, , \,   \epsilon+
|\ul{r}|\, , \, x^{\ul{r}}\, ) \\
  &\beq&  (\, 0\, , \, \dots \, , \, \epsilon -2 \, , \, \epsilon -2 \,
, \,   \epsilon  \, , \, \epsilon \, , \, \dots 
\, , \, \epsilon \, , \,   \epsilon+ |\ul{r}|\, , 
\, x^{\ul{r}}\, ) \, .
  \end{eqnarray*}
  Repeating these steps inductively we arrive at the following form:
  $$
  (\, \epsilon -2  \, , \, \dots \, , \, \epsilon 
-2 \, , \,  \epsilon \, , \, \dots \,
, \, \epsilon \, , \,   \epsilon+ |\ul{r}|\, , \, x^{\ul{r}}\, ) \, \aeq \,
  (\, 1  \, , \, \dots \, , \, 1 \, , \,  3 \, , 
\, \dots \, , \, 3 \, , \,   3+ |\ul{r}|\, , \, 
x^{\ul{r}}\, )\, .
  $$
  If $N=2$ this completes the proof. Otherwise
  we need to transform each pair of the form 
$(x^3y\, , \, x^3y)$ into $(xy\, , \, xy)$.
  Notice that, since $(x^3y)^2=1$, we can move 
this pair everywhere inside the vector
  without changing the other elements. Moreover we 
can conjugate both elements by any of the others
  obtaining again a pair of the form $(g,g)$ with 
$g^2=1$. It follows that we can transform 
$(x^3y\, , \, x^3y)$
  into $(hx^3yh^{-1}\, , \, hx^3yh^{-1})$, for any 
$h\in \langle xy,x^3y,x^{\ul{r}}\rangle=\Dn$,
  hence the result follows (notice that this 
argument follows the proof of Lemma 1.9 in 
\cite{Ka1}).\\

  We consider now the case where $h$ is equal to 
the number of $0$'s in \eqref{dm}.
  In this situation $j$ must be odd, therefore 
there exists an $\ell$ such that 
$GCD(\epsilon:=j+2\ell M,n)=1$,
  where $M=GCD(n,\ul{r})$. From Lemmas \ref{abm} and \ref{pair} we have:
  \begin{eqnarray*}
   (\, 0 \, , \, \dots \, , \, 0 \, , \, j   \, , 
\, \dots \, , \, \, j  \, , \, \, j + |\ul{r}|\,
, \, x^{\ul{r}}\, )&\beq&(\, 0\, , \, \dots \, , \, 0 \, , \, \epsilon \,
, \, \dots \, , \, \epsilon \, , \,   \epsilon+
|\ul{r}|\, , \, x^{\ul{r}}\, ) \\
   &\beq& (\, 0\, , \, \dots \, , \, 0 \, , \, \ell \epsilon \,
  , \,  (\ell +1) \epsilon \, , \, \epsilon \, , \, \dots \, , \, \epsilon \,
, \,   \epsilon+ |\ul{r}|\, , \, x^{\ul{r}}\, ) \\
   \mbox{(for $\ell \epsilon = -1$)} &=&  (\, 0\, , \, \dots \,
  , \, 0 \, , \, -1 \, , \,   \epsilon -1 \, , \, 
\epsilon \, , \, \dots \, , \, \epsilon \, , \, 
\epsilon+ |\ul{r}|\, , \,
x^{\ul{r}}\, ) \\
  &\beq& (\, 0\, , \, \dots \, , \, -\ell \, , \, -\ell -1 \,
, \,   \epsilon -1 \, , \, \epsilon \, , \, \dots 
\, , \, \epsilon \, , \,   \epsilon+ |\ul{r}|\, , 
\, x^{\ul{r}}\, ) \\
  \mbox{(for $\ell = -\epsilon $)} &=& (\, 0\, , \, \dots \,
  , \, \epsilon \, , \, \epsilon -1 \, , \,   \epsilon -1 \,
  , \, \epsilon \, , \, \dots \, , \, \epsilon \, , \,   \epsilon+
|\ul{r}|\, , \, x^{\ul{r}}\, ) \\
  &\beq&  (\, 0\, , \, \dots \, , \, \epsilon -1 
\, , \, \epsilon -1 \, , \,   \epsilon  \,
, \, \epsilon \, , \, \dots \, , \, \epsilon \, , 
\,   \epsilon+ |\ul{r}|\, , \, x^{\ul{r}}\, ) \,.
  \end{eqnarray*}
  Repeating this argument inductively we reach the following form:
  $$
   (\, \epsilon -1\, , \, \dots \, , \, \epsilon -1 \, , \,   \epsilon  \,
  , \, \dots \, , \, \epsilon \, , \,   \epsilon+ 
|\ul{r}|\, , \, x^{\ul{r}}\, ) \, \aeq \,
   (\, 0\, , \, \dots \, , \,0\, , \,   1  \, , \, 
\dots \, , \,1 \, , \,   1+ |\ul{r}|\, , \, 
x^{\ul{r}}\, )\, ,
  $$
  hence the claim follows.\\
  If $h$ is equal to the number of $x^jy$'s, we apply the automorphism
$x\mapsto x$, $x^jy\mapsto y$ and use the equivalence
  $$
(\, -j\, , \, \dots \, , \,-j\, , \,   0  \, , \, \dots \, , \,0 \,
  , \,    |\ul{r}|\, , \, x^{\ul{r}}\, ) \, \beq 
\, (\, 0\, , \, \dots \, , \,0\, , \,   -j  \,
  , \, \dots \, , \,-j \, , \,    |\ul{r}|\, ,
\, x^{\ul{r}}\, )\, .
  $$
  The claim  follows now from the previous case. 
This completes the proof of the Theorem. \qed


\end{document}